\documentclass[12pt,a4paper,twoside]{amsart}
\usepackage{amssymb}
\usepackage{amscd}
\usepackage{eucal}
\usepackage{graphicx}
\usepackage{color}

\title{Logarithmically complete monotonicity of reciprocal arctan function}
\theoremstyle{plain}
 \newtheorem{thm}{Theorem}[section]
 \newtheorem{prop}{Proposition}[section]

\theoremstyle{definition}

\numberwithin{equation}{section}
\author[V.\ Jovanovi\'c]{Vladimir Jovanovi\'c$^1$}
\author[M.\ Treml]{Milanka Treml$^1$}
\address{$^1$Faculty of Sciences and Mathematics,
\newline \indent University of Banja Luka,
\newline \indent Mladena Stojanovi\'ca 2
\newline \indent Banja Luka
\newline\indent Republic of Srpska
\newline\indent Bosnia and Herzegovina}
\email{vladimir.jovanovic@pmf.unibl.org}
\email{milanka.treml@pmf.unibl.org}

\begin{document}
\begin{abstract}
We prove the conjecture stated in F.\ Qi and R.\ Agarwal, \textit{On complete monotonicity for several
classes of functions related to ratios of gamma functions}, J.\, Inequal.\, Appl.\,(2019), 1-42, that the function $1/\arctan$ is logarithmically completely monotonic on $(0,\infty)$, but not a Stieltjes transform.
\end{abstract}

\maketitle


\section{Introduction}

By a \textit{completely monotonic function} (shortly CM) we mean here an infinitely differentiable function $f:(0,\infty)\rightarrow\mathbb{R}$, such that
$$(-1)^nf^{(n)}\ge0,\quad n=0,1,2,\dots.$$
If $f'$ is completely monotonic and $f\ge0$, then we call $f$ a \textit{Bernstein function}. Here we are mostly interested in \textit{logarithmically completely monotonic functions}, that is, infinitely differentiable functions $f:(0,\infty)\rightarrow (0,\infty)$ with the property
\begin{equation*}
(-1)^n(\log f)^{(n)}\ge0,\quad n=1,2,3\dots.
\end{equation*}
A basic fact concerning CM - functions is the Bernstein theorem: a function $f$ is CM if and only if there exists a non-decreasing function $\alpha$ on $(0,\infty)$ satisfying
\begin{equation*}
f(x)=\int_0^{\infty}e^{-xt}d\alpha(t),
\end{equation*}
for all $x>0$ (see \cite{Widder}, p. 161). In some occasions it has been proven a stronger property which leads to complete monotonicity of a function $f$, namely that there exist $a\ge0$ and a non-negative Borel measure $\mu$ on $[0,\infty)$ for which the equality
\begin{equation*}
f(x)=a+\int_0^{\infty}\displaystyle\frac{d\mu(t)}{x+t}
\end{equation*}
holds for $x>0$, where the measure $\mu$ fulfills  the condition
\begin{equation*}
\int_0^{\infty} \frac{d\mu(t)}{1+t}<\infty.
\end{equation*}
Such functions are called \textit{Stieltjes transforms}. We recall that all Stieltjes transforms are logarithmically completely monotonic (see \cite{Berg}), and the latter are CM (see \cite{QiChen}).\\
In \cite{QiAgar} the authors set the conjecture that the function $f(x)=\frac1{\arctan x}$ is logarithmically completely monotonic on $(0,\infty)$, but not a Stieltjes transform. The aim of this paper is to justify these assertions. We will do it in the next section.
\section{Formulations and proofs}
\begin{thm}\label{t1}
The function $f(x)=\frac1{\arctan x}$ is logarithmically completely monotonic on $(0,\infty)$.
\end{thm}
The idea of the proof of Theorem \ref{t1} is based on the Remark (1) in \cite{AlzB}, where the authors suggest employing the residue theorem in an attempt to obtain integral representations of functions under consideration.
\begin{proof}
It suffices to prove that
\begin{equation*}
g(x)=-(\log f(x))'=\displaystyle\frac{1}{(x^2+1)\arctan x}
\end{equation*}
is CM on $(0,\infty)$. In what follows we always assume that $\log$ denotes the principle value of logarithm, i.\ e.\ $\log z=\ln|z|+i\arg z$, with $\arg z\in(-\pi, \pi]$.

\begin{figure}
	\includegraphics[scale=0.3]{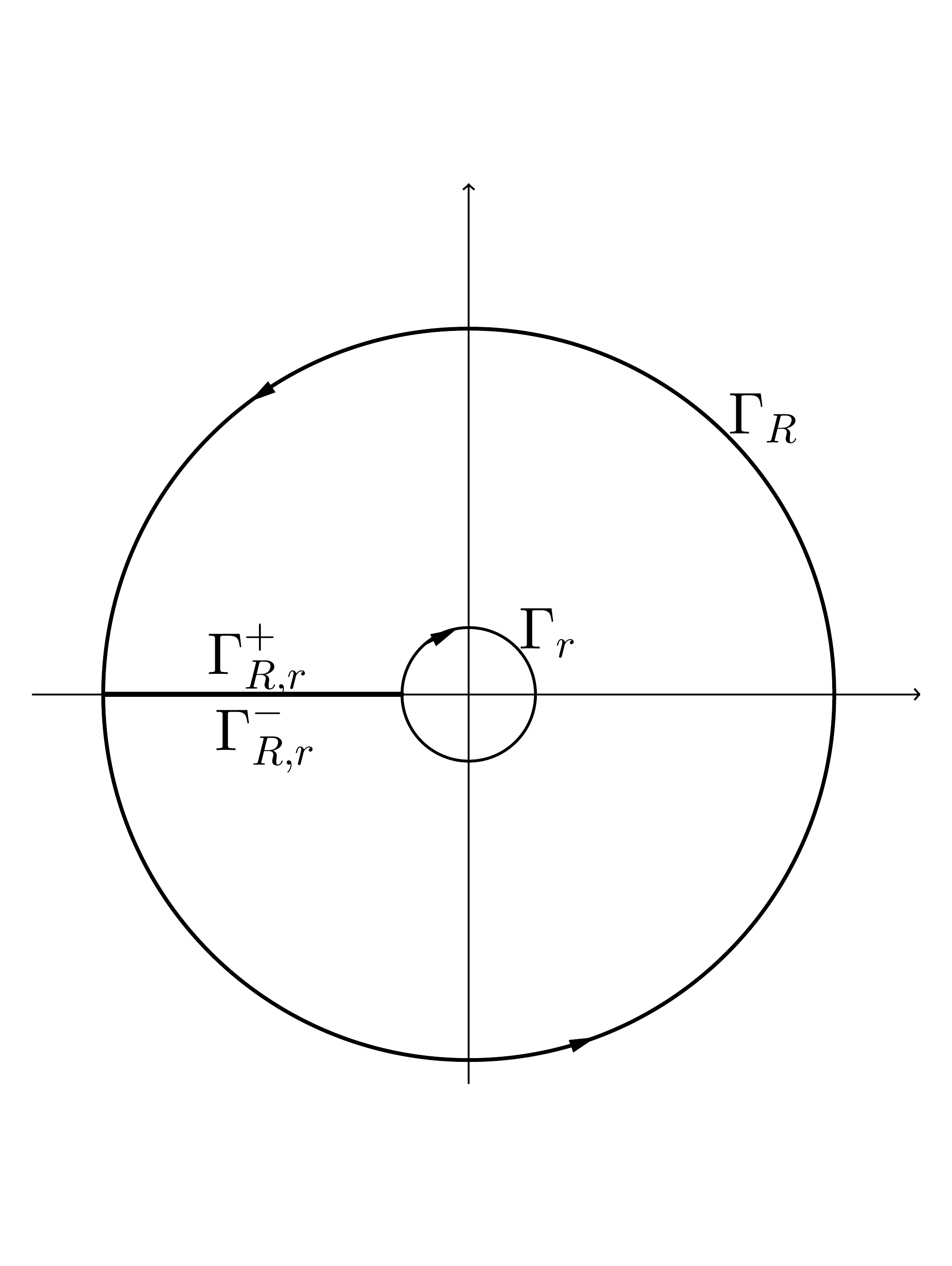}
	\caption{Keyhole contour $\Gamma_{R,r}$}
	\label{kontura}
\end{figure}

Let us consider the integral
$\displaystyle \int_{\Gamma_{R,r}} G(z)\, dz$, over the "keyhole" contour $\Gamma_{R,r}$ given in Figure \ref{kontura}, where
\begin{equation*}
G(z)=\displaystyle\frac{z+1}{z(z-z_0)\log z}
\end{equation*}
and $z_0=\frac{i-x}{i+x}$, for $x>0$. 

We assume $R>1$ and $r<1$. Note that $|z_0|=1$ and that $1,\, z_0$ are the only singularities of $G$ lying inside $\Gamma_{R,r}$. From the residue theorem, we have
\begin{equation*}
\int_{\Gamma_{R,r}}G(z)\, dz=2\pi i(\mathrm{Res}(G(z);z_0)+\mathrm{Res}(G(z);1)).
\end{equation*}
Since $z_0$ is a first-order pole, it follows
\begin{align*}
\mathrm{Res} (G(z);z_0)&=\frac{1+z_0}{z_0\log z_0}=\frac{1+\frac{i-x}{i+x}}{\frac{i-x}{i+x}\log\frac{i-x}{i+x}}=\frac{2i}{(i-x)2i\arctan x}=\\
&=-\frac{(i+x)}{(x^2+1)\arctan x},
\end{align*}
where we used the fact that $\arctan x=\frac{1}{2i}\log \frac{1+ix}{1-ix}$, for $x>0$. Similarly,
\begin{equation*}
\mathrm{Res}(G(z);1)=\lim_{z\to 1}(z-1)\frac{1+z}{z\log z(z-z_0)}=\frac{2}{1-z_0}=\frac{2}{1-\frac{i-x}{i+x}}=\frac{i+x}{x},
\end{equation*}
whence,
\begin{equation}\label{eq: g}
g(x)=\displaystyle\frac{1}{x}-\frac{1}{2\pi i (x+i)}\int_{\Gamma_{R,r}}G(z)\, dz.
\end{equation}
Now, it remains to calculate the integral $\int_{\Gamma_{R,r}}G(z)\, dz$. In order to accomplish it, we start from the relation
\begin{equation}\label{eq: GammaG}
\int_{\Gamma_{R,r}}G(z)\, dz=\int_{\Gamma_R}G(z)\, dz+
\int_{\Gamma_r}G(z)\, dz+\int_{\Gamma_{R,r}^+}G(z)\, dz+\int_{\Gamma_{R,r}^-}G(z)\, dz.
\end{equation}
The first two integrals vanish as $R\to\infty$ and $r\to 0+$. It follows from the estimates
\begin{equation*}
\left|\int_{\Gamma_R}G(z)\, dz\right|\leq 2R\pi\max_{|z|=R}\frac{|z+1|}{|z||\log z||z-x_0|}\leq 2\pi\frac{R+1}{(\ln R-2\pi)(R-1)}
\end{equation*}
and
\begin{equation*}
\left|\int_{\Gamma_r}G(z)\, dz\right|\leq 2r\pi\max_{|z|=r}\frac{|z+1|}{|z||\log z||z-x_0|}\leq
 2\pi\frac{1+r}{(-\ln r-2\pi)(1-r)}.
\end{equation*}
We also have for $t<0$,
\begin{equation*}
\displaystyle \lim_{\genfrac{}{}{0pt}{2}{z\to t}{\Im z>0} 
}G(z)= \frac{t+1}{t(\ln (-t)+\pi i)(t-x_0)}=G^+(t)
\end{equation*}
and
\begin{equation*}
\displaystyle \lim_{\genfrac{}{}{0pt}{2}{z\to t}{\Im z<0} 
}G(z)= \frac{t+1}{t(\ln (-t)-\pi i)(t-x_0)}=G^-(t).
\end{equation*}
Consequently,
\begin{equation}\label{eq: Gamma+}
\int_{\Gamma_{R,r}^+}G(z)\, dz+\int_{\Gamma_{R,r}^-}G(z)\, dz=
\int_{-R}^{-r}[G^-(t)-G^+(t)]\, dt
\end{equation}
Let us denote $\displaystyle I=\lim_{\genfrac{}{}{0pt}{2}{R\to\infty}{r\to 0+}}\int_{\Gamma_{R,r}}G(z)\, dz$. From (\ref{eq: GammaG}) and (\ref{eq: Gamma+}) we obtain
\begin{align*}
I &=\int_{-\infty}^0[G^-(t)-G^+(t)]\, dt\\
  &= \int_{-\infty}^0\frac{-2\pi i(t+1)\, dt}{t(\log^2(-t)+\pi^2)(t-z_0)}\\
  &= -2\pi i\int_{0}^\infty\frac{(1-t)\, dt}{t(\log^2 t+\pi^2)(t+z_0)}.
\end{align*}
Using $z_0=\frac{i-x}{i+x}$, we have
\begin{align*}
 I &= \int_0^{\infty}\frac{-2\pi i(1-t)\, dt}{t(\log ^2 t+\pi^2)(t+\frac{i-x}{i+x})}\\
&=\int_{0}^\infty\frac{-2\pi i(i+x)(1-t)\, dt}{t(\log^2t+\pi^2)(x(t-1)+i(t+1))}\\
&= 2\pi i(i+x)\int_{0}^\infty\frac{((1-t)^2x+i(1-t^2))\,  dt}{t(x^2(1-t)^2+(1+t)^2)(\log^2 t+\pi^2)}.
\end{align*}
Note that (\ref{eq: g}) implies
\begin{equation}\label{eq: gI}
g(x)=\displaystyle\frac{1}{x}-\frac{1}{2\pi i (x+i)}\, I
\end{equation}
and since $\displaystyle\frac{1}{2\pi i (x+i)} I$ is real, we conclude that
\begin{equation*}
\int_{0}^\infty\frac{(1-t^2)\,  dt}{t(x^2(1-t)^2+(1+t)^2)(\log^2 t+\pi^2)}=0.
\end{equation*}
Therefore, from (\ref{eq: gI}), it follows
\begin{equation}
g(x)=\frac{1}{x}-\int_{0}^\infty\frac{(1-t)^2x\,  dt}{t(x^2(1-t)^2+(1+t)^2)(\log^2 t+\pi^2)}.
\end{equation}
Employing
\begin{equation*}
\displaystyle \frac{1}{x}=\int_0^{\infty}\frac{\, dt}{xt(\log^2t+\pi^2)},
\end{equation*}
we get
\begin{equation*}
g(x)=\int_{0}^\infty\frac{(1+t)^2\,  dt}{xt(x^2(1-t)^2+(1+t)^2)(\log^2 t+\pi^2)}.
\end{equation*}
The substitution $t\mapsto\frac1t$ implies
\begin{align*}
\int_{0}^1\frac{(1+t)^2\,  dt}{xt(x^2(1-t)^2+(1+t)^2)(\log^2 t+\pi^2)}=&\\
=\int_{1}^\infty\frac{(1+t)^2\,  dt}{xt(x^2(1-t)^2+(1+t)^2)(\log^2 t+\pi^2)}&.
\end{align*}
Hence,
\begin{equation}\label{eq: g2}
g(x)=2\int_{0}^1\frac{(1+t)^2\,  dt}{xt(x^2(1-t)^2+(1+t)^2)(\log^2 t+\pi^2)}.
\end{equation}
For $a,b,x>0$, it is
\begin{equation*}
\frac{1}{x(a^2x^2+b^2)}=\frac{1}{b^2}\left(\frac1x-\frac12\left(\frac{1}{x+\frac{bi}{a}}+\frac{1}{x-\frac{bi}{a}}\right)\right),
\end{equation*}
and using
\begin{equation*}
\frac{1}{x+\frac{bi}{a}}=
\int_0^\infty e^{-xs}e^{-\frac{bi}{a}s}\, ds, \quad \frac{1}{x-\frac{bi}{a}}=
\int_0^\infty e^{-xs}e^{\frac{bi}{a}s}\, ds
\end{equation*}
and
\begin{equation*}
\frac{1}{x}=\int_{0}^\infty e^{-xs}\, ds,
\end{equation*}
one obtains
\begin{equation*}
\frac{1}{x(a^2x^2+b^2)}=\int_0^{\infty}e^{-xs}\left(1-\cos\frac{bs}{a}\right)\, ds.
\end{equation*}
Setting $a=1-t$ and $b=1+t$ yields
\begin{equation*}
\frac{1}{x(x^2(1-t)^2+(1+t)^2)}=\frac{1}{(1+t)^2}\int_0^\infty e^{-xs}\left(1-\cos\frac{1+t}{1-t}s\right)\, ds. 
\end{equation*}
From (\ref{eq: g2}), we have
\begin{equation*}
g(x)=2\int_0^1\left(\int_0^\infty\frac{e^{-xs}\left(1-\cos\frac{1+t}{1-t}s\right)\, ds}{t(\ln^2 t+\pi^2)}\right)\, dt,
\end{equation*}
and, finally, after interchanging integration order, we obtain
\begin{equation}\label{eq:gfin}
g(x)=\int_0^\infty \left(\int_0^1\frac{2(1-\cos \frac{1+t}{1-t}s)\, dt}{t(\ln^2t +\pi^2)}\right)e^{-xs}\, ds.
\end{equation}
Now, it is evident that (\ref{eq:gfin}) implies complete monotonicity of $g$.
\end{proof}
\begin{thm}\label{t2}
The function $f(x)=\frac1{\arctan x}$ is not a Stieltjes transform on $(0,\infty)$.
\end{thm}
For the proof of this theorem, we the use the following result on Stieltjes transforms from \cite{BergPed}.
\begin{prop}\label{pro}
If $f\neq0$ is a Stieltjes transform, then $1/f$ is a Bernstein function.
\end{prop}
\noindent\textit{Proof of Theorem \ref{t2}}.\\[1ex]
The function $h(x)=\frac1{f(x)}=\arctan x$ is not a Bernstein function, since
\begin{equation*}
h^{(3)}(x)=\displaystyle2\,\frac{1-3x^2}{(1+x^2)^3}
\end{equation*}
changes its sign on $(0,\infty)$. Therefore, according to Proposition \ref{pro}, $f$ if not a Stieltjes transform.\enspace $\Box$


\begin{thebibliography}{1}
%
\bibitem{AlzB}
H.~Alzer and C.~Berg, \textit{Some classes of completely monotonic functions},
  Annales Academiae Scientiarum Fennicae.\ Mathematica \textbf{27} (2002),
  445--460.

\bibitem{Berg}
C.~Berg, \textit{Integral representation of some functions related to the gamma
  function}, Mediterr.\ J.\ Math. \textbf{1} (2004), 433--439.

\bibitem{BergPed}
C.~Berg and H.~L. Pedersen, \textit{A completely monotone function related to
  the gamma function}, Journal of Computational and Applied Mathematics
  \textbf{133} (2001), 219--230.

\bibitem{QiAgar}
F.~Qi and R.~Agarwal, \textit{On complete monotonicity for several classes of
  functions related to ratios of gamma functions}, J Inequal Appl  (2019),
  1--42.

\bibitem{QiChen}
F.~Qi and C.-P. Chen, \textit{A complete monotonicity property of the gamma
  function}, J.\ Math.\ Anal.\ Appl. \textbf{296} (2004), 603--607.

\bibitem{Widder}
D.~V. Widder, \textit{The Laplace transform}, Princeton University Press, 1946.
%
\end{thebibliography}
\end{document}